\documentclass[a4paper]{amsart}
\usepackage[utf8]{inputenc}
\usepackage[T1]{fontenc}

\usepackage[style=alphabetic,maxbibnames=10,backend=bibtex]{biblatex}
\appto{\bibsetup}{\emergencystretch=0.75em}

\usepackage{hyperref}

\makeatletter
\newcommand{\proofstep}[1]{%
  \par%
  \addvspace{\medskipamount}%
  \textit{#1\@addpunct{.}}\enspace\ignorespaces
}  
\makeatother
	
\newenvironment{acknowledgement}%
  {\par%
   \hspace{\parindent}\textsc{Acknowledgements.}\enskip\ignorespaces}
  {\par}

\usepackage{amsmath, amssymb, amsfonts, mathtools, amsthm, mathrsfs, dsfont, fge, extpfeil}
\usepackage{tikz,tikz-cd}
\usetikzlibrary{matrix,arrows,positioning,calc}

\newcommand\Q{\mathbb{Q}} 
\newcommand\Z{\mathbb{Z}}
\newcommand\N{\mathbb{N}}
\newcommand\M{\mathcal{M}}
\newcommand\cN{\mathcal{N}}
\newcommand\T{\mathbb{T}}
\newcommand\bS{\mathbb{S}}
\newcommand\del{\partial}
\newcommand\J{\mathit{J}}

\newcommand\C{\mathcal{C}}
\newcommand\cR{\mathcal{R}}
\newcommand\cP{\mathcal{P}}
\newcommand\cQ{\mathcal{Q}}

\newcommand\id{\mathrm{id}}
\DeclareMathOperator{\rad}{rad}

\DeclareMathOperator{\im}{im}

\def\setminus{\mathbin{\fgebackslash}}
\DeclareMathOperator{\rep}{rep}

\DeclareMathOperator{\mods}{mod}

\DeclareMathOperator{\ind}{ind}
\DeclareMathOperator{\add}{add}
\DeclareMathOperator{\Hom}{\mathrm{Hom}}
\DeclareMathOperator{\End}{\mathrm{End}}
\DeclarePairedDelimiter{\card}{|}{|}

\DeclarePairedDelimiter{\ceil}{\lceil}{\rceil}
\DeclareMathOperator{\dimv}{\underline{\dim}}

\def\Ext#1#2{\mathrm{Ext}_{#1}^{#2}}
\def\defined#1{\textbf{#1}}
\newcommand{\iso}{\xrightarrow{\sim}}
\def\isom{\mathrel{\cong}}

\theoremstyle{definition}
\newtheorem{mydef}{Definition} 

\theoremstyle{plain}
\newtheorem{theorem}[mydef]{Theorem}
\newtheorem*{theorem*}{Theorem}
\newtheorem*{maintheorem}{Main Theorem}
\newtheorem{prop}[mydef]{Proposition}
\newtheorem{proposition}[mydef]{Proposition}
\newtheorem{lemma}[mydef]{Lemma}

\newtheorem*{corollary*}{Corollary}

\theoremstyle{remark}
\newtheorem{example}[mydef]{Example}
\newtheorem{examples}[mydef]{Examples}
\newtheorem*{example*}{Example}

\newtheorem*{remark}{Remark}

\newtheorem*{remarks}{Remarks}

\title{Tame hereditary path algebras and amenability}
\author{Sebastian Eckert}
\address{Sebastian Eckert, Fakultät für Mathematik, Universität Bielefeld, Postfach 100~131, 33501~Bielefeld, Germany.}
\email{seckert@math.uni-bielefeld.de}
\thanks{The author has been supported~by the Alexander von~Humboldt Foundation in the framework of an Alexander von Humboldt Professorship endowed by the German Federal Ministry of Education and Research.}


\begin{document}

\begin{abstract}
In this note we are concerned with the notion of amenable representation type as defined in a recent paper by Gábor Elek. We will show that the tame hereditary path algebras of quivers of extended Dynkin type over any field~$k$ are of amenable type, thus extending a conjecture in the aforementioned paper to another class of tame algebras. In doing so, we avoid using already known results for string algebras. We also show that path algebras of wild acyclic quivers over finite fields are not amenable.
\end{abstract}

\maketitle

\section{Introduction}

In \cite{Elek2017InfiniteDimensionalRepresentationsAmenabilty}, Elek has introduced the notion of hyperfiniteness for countable sets of modules and that of amenable representation type for algebras. We deviate slightly and give the following definiton, as we extend the notion of hyperfiniteness to arbitrary families and that of amenability to representation-finite algebras.

\begin{mydef} \label{def:HyperfinitenessAmenability}
Let~$k$ be a field, $A$ be a finite dimensional $k$-algebra and let $\M$ be a set of $A$-modules. $\M$ is called \defined{hyperfinite} provided for every $\varepsilon > 0$ there exists a number $L_\varepsilon > 0$ such that for every $M \in \M$ there exists a submodule $P \subseteq M$ such that
\begin{equation} \label{eq:HFSubmoduleBig} \dim_k P \geq (1-\varepsilon) \dim_k M, \end{equation} and modules $N_1,N_2, \dots N_t \in \mods A$, with $\dim_k N_i \leq L_\varepsilon$, such that $P \isom \bigoplus_{i=1}^{t} N_i$.

The $k$-algebra $A$ is said to be of \defined{amenable representation type} provided the set of all finite dimensional $A$-modules (or more specific, a set which meets any isomorphism class of finite dimensional $A$-modules) is hyperfinite.
\end{mydef}

Using this definition, families of modules are classified by the existence of large submodules which are built up from modules of bounded dimension. By doing so, one can divide those alegbras into amenable and non-amenable representation types, possibly mirroring the tame/wild dichotomy. It is further conjectured \cite[Conjecture~1]{Elek2017InfiniteDimensionalRepresentationsAmenabilty} that finite dimensional algebras are of tame type if and only if they are of amenable representation type. While string algebras (over countable fields) are shown to be of amenable type in \cite[Proposition~10.1]{Elek2017InfiniteDimensionalRepresentationsAmenabilty}, further examples are still lacking.

In this note we prove the conjecture for path algebras of acyclic quivers, exhibiting a new class of algebras of tame and amenable representation type.
\begin{maintheorem}
Let~$Q$ be an acyclic quiver of extended Dynkin type $\widetilde{A}_n$, $\widetilde{D}_n$, $\widetilde{E}_6$, $\widetilde{E}_7$ or $\widetilde{E}_8$. Let~$k$ be any field. Then the path algebra~$kQ$ of~$Q$ is of amenable representation type.

Moreover, if $k$ is a finite field and $Q$ is an acyclic quiver of wild type, then the path algebra~$kQ$ of~$Q$ is not of amenable type.
\end{maintheorem}

In order to prove our Main Theorem, we will first deduce some general results for hyperfinite families in Section~\ref{section:Hyperfiniteness}.
In Section~\ref{section:LemmasExtendedDynkin}, we will then take the first steps towards proving the Main Theorem by discussing several technical lemmas, before studying the case of the $2$-Kronecker quiver in Section~\ref{section:KroneckerQuiver}. This result will then be used to prove the first half of our Main Theorem in Section~\ref{section:ExtendedDynkinAmenable}, using a descent argument via localization.
The other direction will be shown in Section~\ref{section:WildQuivers}, using a similar simplification and suitable functors.

\section{Hyperfiniteness and Amenability} \label{section:Hyperfiniteness}
We shall start the discussion by a further inspection of Definition~\ref{def:HyperfinitenessAmenability}.

\begin{remark}
Since $P$ is a submodule of~$M$, the condition in (\ref{eq:HFSubmoduleBig}) in the above definition is equivalent to
\begin{equation} \label{eq:HFQuotientSmall} \dim_k (M / P) \leq \varepsilon \dim_k M, \end{equation}
for $\dim_k (M / P) = \dim_k M - \dim_k P$.
\end{remark}

\begin{remark}
Finite sets are hyperfinite, as we can take $L_\varepsilon$ to be the maximum of the dimensions. For the same reason, families of modules of bounded dimension are hyperfinite.
Moreover, if $\M, \M'$ are hyperfinite families, so is $\M \cup \M'$: we can simply choose $L_\varepsilon$ to be the maximum of the $L_\varepsilon$s, corresponding to $\M$ and to $\M'$ respectively. Similarly, any finite union of hyperfinite families is hyperfinite.
\end{remark}

\begin{proposition} \label{prop:AdditiveClosureStaysHyperfinite}
Let $\M$ be a family of $A$-modules. If $\M$ is hyperfinite, so is the family of all finite direct sums of modules in $\M$.
\begin{proof}
Let $\M$ be hyperfinite and let $X = \bigoplus_{i=1}^{n} M_i$.
Now let $\varepsilon > 0$. Choose $L_\varepsilon$ as one would to show the hyperfiniteness of~$\M$.
For each $1 \leq i \leq n$, choose \[\bigoplus_{j=1}^{t_i} N_{i,j} = N_i \subseteq M_i,\] as for the hyperfiniteness of~$\M$.
Then \[N := \bigoplus_{i = 1}^{n} N_i \subseteq \bigoplus_{i = 1}^{n} M_i = X,\] as direct sums respect submodule inclusions.
Also \[\dim_k N = \sum_{i = 1}^{n} \dim_k N_i \geq \sum_{i = 1}^{n} (1-\varepsilon) \dim_k M_i = (1-\varepsilon) \sum_{i = 1}^{n} \dim_k M_i = (1-\varepsilon) \dim_k X.\]
Moreover, $\dim_k N_{i,j} \leq L_\varepsilon$.
\end{proof}
\end{proposition}

This shows that to check amenability, it is enough to check the criterion on all indecomposable modules.

\begin{example}
Let $A$ be a representation finite algebra. Since there are only finitely many isoclasses of (finitely generated) indecomposable modules, for any $\varepsilon$ by choosing
\[ L_\varepsilon = \max\limits_{M \in \ind A} \{\dim_k M\},\] for any $M \in \ind~A$ we can just take $P=M$ to show that such an algebra $A$ is of amenable representation type.\\
A non-example is given in \cite[Theorem~6]{Elek2017InfiniteDimensionalRepresentationsAmenabilty} by the wild Kronecker algebras over countable fields, while the 2-Kronecker algebra is of amenable representation type (see also Theorem \ref{thm:2KroneckerAmenable}).
\end{example}

\begin{proposition} \label{prop:ExtendingHFfromSubmodulesOfBoundedCodimension}
Let $A$ be a finite dimensional $k$-algebra. Let $\M,\cN$ be two sets of finite dimensional $A$-modules, where $\cN$ is hyperfinite. If there is some $L \geq 0$ such that for all $M \in \M$, there exists a submodule $P \subset M$ of codimension less than or equal to~$L$ with $P \in \cN$, then $\M$ is also hyperfinite.
\begin{proof}
Let $L \geq 0$ and $\varepsilon > 0$. If the dimension of the modules in $\M$ was bounded, say by~$K$, we can set $L_\varepsilon := K$ and choose $P = M$  for all $M \in \M$. Hence, let $M \in \M$ have large dimension, i.e. $\dim_k M > \frac{2L}{\varepsilon}$. We choose a submodule $P \in \cN$ of codimension bounded by $L$. Since $\cN$ is hyperfinite, there is some submodule $Y \subseteq P$ such that $\dim Y \geq (1-\frac{\varepsilon}{2})\dim P$, while $Y$ decomposes into direct summands of dimension less than or equal to~$L'_{\frac{\varepsilon}{2}}$.

We thus have that
\begin{align*}
\dim Y &\geq \left(1-\frac{\varepsilon}{2}\right) \dim P = \dim P - \frac{\varepsilon}{2} \dim P\\
&\geq \left( \dim M - L \right) - \frac{\varepsilon}{2} \dim M \\
&\geq \dim M - \frac{\varepsilon}{2} \dim M - \frac{\varepsilon}{2} \dim M\\
&= (1-\varepsilon) \dim M,
\end{align*}
using the fact that for~$M$ we have that $L \leq \frac{\varepsilon}{2} \dim M$.
What is more, $Y$~decomposes into direct summands of dimension less than or equal to~$L'_{\frac{\varepsilon}{2}}$. If we therefore choose $L_\varepsilon$ to be the maximum of $L'_{\frac{\varepsilon}{2}}$ and $\frac{2L}{\varepsilon}$, we have shown that $\M$~is hyperfinite.
\end{proof}
\end{proposition}

\begin{proposition} \label{prop:HFPreservingFunctors}
Let~$k$ be a field and $A,B$~be two finite dimensional $k$-algebras. Let $F \colon \mods A \to \mods B$ be an additive, left-exact functor such that there exists $K_1, K_2 > 0$~with
\begin{equation} \label{eq:EquivalenceCondHFPresFun}
K_1 \dim X \leq \dim F(X) \leq K_2 \dim X,\end{equation}
 for all $X \in \mods A$.
If $\cN \subseteq \mods A$ is a hyperfinite family, then the family \[\{F(X) \colon X \in \cN\} \subseteq \mods B\] is also hyperfinite.
\begin{proof}
Let us denote $\M := \{F(X) \colon X \in \cN\}$.
By the hypothesis, for any $\tilde{\varepsilon}$ we can find some $\tilde{L}_{\tilde{\varepsilon}} > 0$ to exhibit the hyperfiniteness of the family $\cN$.
Let $M \in \mods B$ such that there exists some $N \in \cN$ with $F(N) = M$. Then there is a submodule $P \subseteq N$ such that
$P = \bigoplus_{i=1}^{t} P_i$ with $\dim P_i \leq \tilde{L}_{\tilde{\varepsilon}}$ and $\dim P \geq (1-\tilde{\varepsilon}) \dim N$.
Since $F$ is additive, we have that $F(P) = \bigoplus_{i=1}^{t} F(P_i)$, and \[\dim F(P_i) \leq K_2 \dim P_i \leq K_2 \tilde{L}_{\tilde{\varepsilon}}.\]
Moreover, the sequence \[0 \to F(P) \to M \to F(N/P)\] is exact, so $F(P)$ is a submodule of~$M$, and by the rank-nullity theorem,
\begin{align*} \dim F(P) &\geq \dim M - \dim F(N/P) \\
& \geq \dim M - K_2 \dim N/P \\ 
& \geq 
\dim M - K_2 \tilde{\varepsilon} \dim N \\
& \geq \dim M - \frac{K_2}{K_1} \tilde{\varepsilon} \dim M = (1-\varepsilon) \dim M, \end{align*}
if we choose $\tilde{\varepsilon} = \frac{K_1}{K_2} \varepsilon$.
We can therefore choose $L_\varepsilon$ to be $K_2 \tilde{L}_{\tilde{\varepsilon}}$ to show the hyperfiniteness of~$\M$.
\end{proof}
\end{proposition}

\begin{remarks}
A functor fulfilling the hypothesis of Proposition~\ref{prop:HFPreservingFunctors} may be called hyperfiniteness preserving or HF-preserving.
Moreover, inspection of the proof shows that the left inequality of (\ref{eq:EquivalenceCondHFPresFun}) need only hold for $X \in \cN$.
\end{remarks}

\begin{example} The application of Proposition~\ref{prop:HFPreservingFunctors} includes some useful cases:
\begin{enumerate}
\item Equivalences are HF-preserving functors, since they are left exact and simples getting mapped to simples determine length and thus dimension.
\item If $A$ is the path algebra of a quiver~$Q$ of amenable representation type, and $i \in Q_0$ is a sink, the path algebra~$kQ'$, where $Q' = \sigma_i(Q)$ is the quiver obtained from~$Q$ by reversing all arrows starting or ending in $i$, is also of amenable type. To see this, let $\C = \C_i$ be the full subcategory of $\mods kQ'$ of objects having no direct summand isomorphic to $S(i)$. Then every indecomposable object of $\mods kQ'$ is either contained in $\C$ or is isomorphic to $S(i)$. Consider the BGP reflection functor $S_i^{+}$, which is left exact. For $S_i^{+}$ and $S_i^{-}$, the right hand part of (\ref{eq:EquivalenceCondHFPresFun}) holds with $K_2 = \card*{Q_1}+1$, implying that the right inequality holds for $S_i^{+}$ on $\C$. As $\C$ is in the essential image of~$S_i^{+}$, the family of indecomposable modules of $\mods kQ'$ is hyperfinite, proving the claim.
\end{enumerate}
\end{example}

\section{Extended Dynkin quivers} \label{section:LemmasExtendedDynkin}

From now on, let~$k$ be any field and~$Q$ be an acyclic, extended Dynkin quiver, i.e. a quiver of type $\widetilde{A}_n$, $\widetilde{D}_n$, $\widetilde{E}_6$, $\widetilde{E}_7$ or $\widetilde{E}_8$ which has no oriented cycles. Let $Q_0$ be the underlying set of vertices and $Q_1$ the set of arrows, and denote the starting and terminating vertex of an arrow $a \colon i \to j$ by $s(a) = i$ and $t(a) = j$ respectively.
We shall consider the category $\mods kQ$ of finite dimensional right modules over the path algebra~$kQ$, which is equivalent to the category $\rep kQ$ of finite dimensional $k$-linear representations of~$Q$. Here, a representation $(M_i,M_a)$ of~$Q$ is given by a collection $\{M_i \colon i \in Q_0\}$ of finite dimensional $k$-vector spaces along with a collection of $k$-linear maps $\{M_a \colon a \in Q_1\}$. It is well known that these categories are tame and hereditary. We shall write $\dimv M$ for the \defined{dimension vector} of the representation $M$, where $(\dimv M)_i = \dim M_i = \dim_k M_i$.

Recall the \defined{Euler bilinear form},
\[\langle x, y \rangle := \sum_{i \in Q_O} x_i y_i - \sum_{a \in Q_1} x_{s(a)} y_{t(a)},\]
which agrees with the homological bilinear form on the dimension vectors of the representations, i.e.
\[\langle \dimv X, \dimv Y \rangle = \sum_{t\geq 0} \dim_k \Ext{kQ}{t}(X,Y),\]
the \defined{Tits form}, i.e. the corresponding quadratic form, $q(x) = \langle x,x\rangle$, and its radical 
\[\rad q = \{x\in\Q^{\card*{Q_0}} \colon q(x) = 0 \}.\]
We shall denote the minimal (integer) generator of $\rad q$ by~$h_Q$.
What is more, define the \defined{defect} of~$X$ by $\del(X) = \del(\dimv X) = \langle h_Q, \dimv X\rangle$. It allows us to distinguish between preprojective ($\del < 0$), regular ($\del = 0$) and preinjective ($\del > 0$) components of $\mods kQ$.

Also recall that $\mods kQ$ has Auslander-Reiten sequences, giving rise to the AR translate $\tau$ and its inverse $\tau^-$. The category may thus be described by its AR quiver~$\Gamma_{kQ}$. Moreover, there exists a transformation $c$ on $\Q^{\card*{Q_0}}$, called the \defined{Coxeter transformation}, such that for any module without projective direct summand, we have 
\[\dimv \tau(X) = c(\dimv X).\]
Note that the defect is a uniquely determined, normalised $c$-invariant form.
Now, the \defined{Coxeter number} $d_Q$ is the smallest positive integer $d$ such that for all $x \in~\Q^{\card*{Q_0}}$,
\[ c^{d}(x) - x \in \rad q. \]

What is more, we recall that the regular component of~$\Gamma_{kQ}$ is an Abelian category comprising pairwise orthogonal standard stable tubes, of which at most three are inhomogeneous, having more than one simple regular module. As the simple regular modules of each tube form a cycle under $\tau$, we may therefore use a triple $(p,q,r)$ of positive integers to list the cycle lengths of the inhomogeneous tubes and call it the \defined{tubular type} of~$Q$.

We will further use the notion of a \defined{perpendicular category} for some module $M \in \mods kQ$ by
\[X^\perp := \{Y \in \mods~kQ \colon \Hom(X,Y) = 0 = \Ext{}{1}(X,Y)\}.\]

To proceed with the proof of the main theorem, we gather some technical lemmas. The first will be a result on the sum of the dimension vectors of the simple regular modules in homogeneous tubes of~$\Gamma_{kQ}$.

\begin{lemma} \label{lemma:SumRegularSimplesIshQ}
Let~$k$ be any field.
Let $\T$ be an inhomogeneous tube of rank~$m$ in~$\Gamma_{kQ}$. Let $S_1, \dots, S_m$ denote the isoclasses of simple regular modules on the mouth of~$\T$, such that $\tau S_i = S_{i-1}$ for $i = 2, \dots, m$ and $\tau S_1 = S_m$.
Then $\sum_{i=1}^{m} \dimv S_i = h_Q$.
\begin{proof}
For the case of an algebraically closed field, we may argue as follows: By \cite[Theorem~3.6.(5)]{Ringel1984TameAlgebrasIntegralQuadraticForms} in connection with \cite[Section~3.4]{Ringel1984TameAlgebrasIntegralQuadraticForms}, the regular component of~$\Gamma_{kQ}$ is given by certain tubes $\T(\rho)$, which are each generated by orthogonal bricks $E_1^{(\rho)}, \dots, E_{m_\rho}^{(\rho)}$. As these lie on the mouth of a stable tube of rank $m_\rho$, the isoclasses $S_1, \dots, S_m$ correspond to the $E_i^{(\rho)}$ for some $\rho$.
Now it follows from an argument in the proof of \cite[3.4.(10)]{Ringel1984TameAlgebrasIntegralQuadraticForms} that $\sum_{i=1}^{m} \dimv S_i = h_Q$. 

For the general case, one may inspect the relevant tables of \cite[Chapter~6]{DlabRingel1976IndecomposableRepresentationsGraphsAlgebras}, where the dimension vectors of the simple regular modules are listed.
\end{proof}
\end{lemma}

\begin{lemma} \label{lemma:RegularsInTubeRankGreaterOne}
Let $\T$ be a tube of rank $m \geq 2$ in~$\Gamma_{kQ}$. Let $X$ be an indecomposable regular in $\T$.
Then there exists a submodule $Y \subseteq X$ of codimension bounded by the sum of the entries of~$h_Q$ and a simple regular module $T \in \T$ such that $Y \in T^{\perp}$.
\begin{proof}
Let $T_1, \dots, T_m$ denote the isoclasses of simple regular modules on the mouth of~$\T$ such that $\tau T_i = T_{i-1}$ for $i = 2, \dots, m$ and $\tau T_1 = T_m$.
Following \cite[Chapter~3]{Ringel1984TameAlgebrasIntegralQuadraticForms}, we define the objects $T_i[\ell]$. First, let $T_i[1] := T_i$ for each $1 \leq i \leq m$. Now, for $\ell \geq 2$, recursively define 
$T_i[\ell]$ to be the indecomposable module in $\T$ with $T_i[1]$ as a submodule such that $T_i[\ell] / T_i[1] \isom T_{i+1}[\ell-1]$. Thus $T_i[\ell]$ is the regular module of regular length $\ell$ with regular socle $T_i$. Now, any regular indecomposable in $\T$ will be given as some $T_i[\ell]$. We may define $T_i[\ell]$ for all $i \in \Z$ by letting $T_i[\ell] \isom T_j[\ell]$ iff $i \equiv j \mod~m$.
Note that \begin{equation} \label{eq:DimvIndecReg} \dimv T_i[\ell] = \sum_{j=i}^{i+\ell-1} \dimv T_j.\end{equation}

By \cite[3.1.(3')]{Ringel1984TameAlgebrasIntegralQuadraticForms} and its extension to an arbitrary field via \cite[Corollary~1]{Ringel1994BraidGroupActionExceptionalSequences} we have that
\[ \langle \dimv T_i, \dimv T_j \rangle =
\begin{cases}
\phantom{-}1, & i \equiv j \mod~m, \\
-1, & i \equiv j+1 \mod~m, \\
 \phantom{-}0, & \text{else.}
\end{cases} \] 
This implies that for any $j \in \Z$, $\langle \dimv T_j, \dimv T_i[\ell]\rangle = 0$, provided $\ell \equiv 0 \mod~m$. For $T_i[\ell]$ is uniserial, we have that $\Hom(T_j,T_i[\ell]) = 0$ if and only if $j \not\equiv i \mod~m$. Since $m \geq 2$, this implies that for all $i \in \Z$ and $\ell \equiv 0 \mod~m$, $T_i[\ell]$ is contained in the perpendicular category of some simple regular, e.g. in $T_{i+1}^\perp$.

If $\ell \not\equiv 0 \mod~m$, then write $\ell = n \cdot m + r$, where $0 < r < m$. Then there is a short exact sequence
\[ 0 \to T_i[nm] \to T_i[\ell] \to Z \to 0,\]
where $\dimv Z \leq h_Q$ using (\ref{eq:DimvIndecReg}) and Lemma~\ref{lemma:SumRegularSimplesIshQ}. Thus, we have found a suitable submodule $T_i[nm] \in T_{i+1}^\perp$. 
\end{proof}
\end{lemma}

\begin{lemma} \label{lemma:PPofDefectSmallerThanTubeRank} 
Let $X=\tau^{-r} P(i)$ be some indecomposbale preprojective $kQ$-module of defect $\del(X) = -d < 0$. Let $\T$ be a tube of rank $m > d$. Then there is a simple regular $S \in \T$ such that $X \in S^\perp$.
\begin{proof}
Clearly,
\[-d = \langle h_Q, \dimv X\rangle = - \langle \dimv X, h_Q\rangle = -\langle \dimv P(i), h_Q\rangle = -(h_Q)_i,\]
so $h_Q$ has a component equal to~$d$.
Let $S_1, \dots, S_m$ be the simple regular modules on the mouth  of~$\T$. 
By Lemma \ref{lemma:SumRegularSimplesIshQ}, $\sum_{j=1}^{m} \dimv S_j = h_Q$. Now, since $m > d$, only up to $d$ out of the $m$ modules can have a vector space at vertex $i$ which is non-zero. Let $j_0$ be a such that $(\dimv S_{j_0})_i = 0$.
For $1 \leq j \leq m$,
\begin{align*}
\langle \dimv S_j, \dimv \tau^{-r} P(i) \rangle &= -\langle \dimv \tau^{-r} P(i), c(\dimv S_j)\rangle = - \langle \dimv P(i), c^{r+1}(\dimv S_j)\rangle \\
& = \langle \dimv P(i), \dimv S_{j-r-1} \rangle = (\dimv S_{j-r-1})_i.
\end{align*}
Now, we can choose $j$ such that $j-r-1 \equiv j_0 \mod~m$.
Since $\Hom(S_j,X)=0$ , for $S_j$ is regular and $X$ is preprojective, we must have that $\Ext{}{1} (S_j,X) = 0$. Thus, $X \in S_j^{\perp}.$
\end{proof}
\end{lemma}

Indeed, a slightly stronger result can be shown.

\begin{lemma} \label{lemma:PPofDefectSmallerThanTwiceTubeRank}
Let $X=\tau^{-r} P(i)$ be some indecomposbale preprojective $kQ$-module of defect $\del(X) = -d < 0$ and $r>d_Q$.
Let $\T$ be a hon-homogeneous tube of rank $m > \frac{d}{2}$. Then one of the following holds:
\begin{enumerate}
\item There exists a simple regular module $S \in \T$ such that $X \in S^{\perp}$.
\item There exists a submodule $Y \subseteq X$ and simple regular modules $S,T \in \T$ such that $0 \to Y \to X \to T \to 0$ is exact and $Y \in S^{\perp}$.
\end{enumerate}
\begin{proof}
Clearly,
\[-d = \langle h_Q, \dimv X\rangle = - \langle \dimv X, h_Q\rangle = -\langle \dimv P(i), h_Q\rangle = -(h_Q)_i,\]
so $h_Q$ has a component equal to~$d$.
Let $S_1, \dots, S_m$ be the simple regular modules on the mouth of~$\T$.
Then $\sum_{j=1}^{m} \dimv S_j = h_Q$ by Lemma \ref{lemma:SumRegularSimplesIshQ}. We will write $d_j = (\dimv S_{j})_i$ and have $\sum_{j = 1}^{m} d_j = d$.

Now, for $1 \leq j \leq m$, \begin{align*} \dim \Hom(X,S_j) - \underbrace{\dim \Ext{}{1}(X,S_j)}_{=0} &= \langle \dimv X, \dimv S_j \rangle = \langle \dimv \tau^{-r} P(i), \dimv S_j \rangle \\
& = \langle \dimv P(i), \dimv S_{j-r} \rangle = (\dimv S_{j-r})_i = d_{j-r}.\end{align*}
where $\dim \Ext{}{1}(X,S_j) = \dim D \Hom(S_{j+1}, X) = 0$ by \cite[2.4.(6*)]{Ringel1984TameAlgebrasIntegralQuadraticForms} and the fact that there are no maps from the regular to the preprojective component.

Thus, if there is some $j$ such that $(\dimv S_{j-r})_i > 0$, we can choose some non-zero map $\theta \in \Hom(X,S_{j})$. The image $\im \theta \subseteq S_{j}$ must be regular or has a preprojective summand. If there was a preprojective summand $Z$, it must be to the right of~$X$ in~$\Gamma_{kQ}$. But for any preprojective module $M$ in the $r$-th translate of the projectives or further to the right, we know that $\dimv M = \dimv \tau^{d_Q} M - m h_Q > h_Q$, where $m < 0$ by the properties of preprojective modules.
On the other hand, $\dimv \im \theta \leq \dimv S_{j} \leq h_Q$, a contradiction. Thus, $\im \theta$ must be regular. Since $S_{j}$ is a simple regular, this implies that $\theta$ is surjective. We therefore have an exact sequence
\[0  \to Y \to X \to S_{j} \to 0,\] by letting $Y := \ker \theta$.
Applying $\Hom(-,S_{j})$, we get an exact sequence
\[\xi\colon 0 \to \Hom(S_{j},S_{j}) \to \Hom(X,S_{j}) \to \Hom(Y,S_{j}) \to \Ext{}{1}(S_{j},S_{j}) \to \Ext{}{1}(X,S_{j}).\]
Since $S_{j}$ is a non-homogeneous simple regular, there are no self-extensions, and we have $\Ext{}{1}(S_{j},S_{j}) = 0$. Hence, $\xi$ yields the short exact sequence
\[\xi' \colon 0  \to \End(S_j) \to \Hom(X,S_j) \to \Hom(Y,S_j) \to 0.\]

Now, if $d_{j-r} > 1$ for all $j$, using the hypothesis, we would have 
\[2m \leq \sum_{j=1}^{m} d_j = d < 2m,\] a contradiction.
Therefore, there is some $j_0$ with $d_{j_0-r} \leq 1$.

If $d_{j_0-r} = 1$, we have that $\dim \Hom(X,S_{j_0}) = 1$, so the exact seqeunce $\xi'$ implies that $\Hom(Y,S_{j_0})$ must be zero by dimension arguments, for $\dim \End(S_{j_0}) \geq 1$.
Hence, \[\Ext{}{1}(S_{j_0+1},Y) =  D \Hom(Y,S_{j_0}) = 0.\] Along with the fact that $\Hom(S_{j_0+1},Y) = 0$, since there are no maps from regular to preprojective modules, this implies that $Y \in  S_{j_0+1}^{\perp}$ and we are in case (2).

If $d_{j_0-r} = 0$, we have that $\dim \Hom(X,S_{j_0}) = 0$. So, similarly, \[ \Ext{}{1}(S_{j_0+1},X) = 0, \quad \text{and} \quad \Hom(S_{j_0+1},X) = 0,\]
and $X \in S_{j_0+1}^{\perp}$, showing that we are in case (1).
\end{proof}
\end{lemma}

\begin{remark}
Similar results can reasonably expected for preinjective modules. Yet, we choose a different approach for those.
\end{remark}

\begin{lemma} \label{lemma:DescentOnDefectOfPIs}
Let $X = \tau^{r} I(i)$ be some indecomposable preinjective $kQ$-module of defect $\del(X) = d$. Then there is an injective module $I(j)$ such that there exists a non-zero morphism $\theta \colon X \to I(j)$. Moreover,
for any direct summand $Z$ of $\ker \theta$, we have $\del(Z) < d$.
\begin{proof}
Let $E(X)$ be the injective envelope of~$X$, and take $I(j)$ to be some indecomposable direct summand of~$E(X)$. This yields a non-zero homomorphism $\theta \colon X \to E(X) \to I(j)$.
Consider the exact sequence
\[0 \to \ker \theta \to X \to \im \theta \to 0.\]
Since there is a map from a preinjective module to $\im \theta$, the latter must be preinjective or zero. Yet, $\im \theta \neq 0$, since $\theta$ is non-zero. Thus, $\del(\im \theta) > 0$. 
Therefore, $\del(\ker \theta) < \del(X)$.
If $\ker \theta$ only had preprojective or regular summands $Z$, we are done, for then $\del(Z) \leq 0$.
Thus, we may assume that there is some some preinjective direct summand $Z$. Since $Z$ embeds into $\ker \theta$ and the kernel embeds into $X$, we get a short exact sequence
\[0 \to Z \to X \to X/Z \to 0,\]
using the fact that $\mods kQ$ is Abelian.
Since $X$ is preinjective, again $X/Z$ must be preinjective or zero. If it was zero, then $Z \isom \ker \theta \isom X$, a contradiction, since $\del(\ker \theta) \neq \del(X)$. Thus, $\del(X/Z) > 0$, and hence we may conclude
$\del(Z) < \del(X) = d$.
\end{proof}
\end{lemma}

\section{The $2$-Kronecker quiver case} \label{section:KroneckerQuiver}

The next theorem will be used as the base case in the proof of the main result. In the case of a countable algebraically closed field, it follows from the results on string algebras in \cite[Proposition~10.1]{Elek2017InfiniteDimensionalRepresentationsAmenabilty}, but we give a direct and independent proof here for illustration purposes and to the convenience of the reader.

\begin{theorem} \label{thm:2KroneckerAmenable}
Let~$k$ be any field. Then the path algebra of the $2$-Kronecker quiver is of amenable representation type.
\begin{proof}
We fix a notation for the  the vertices and arrows as follows.
\begin{center}
\begin{tikzpicture}
\matrix (m) [matrix of math nodes,row sep=1em,column sep=3em,minimum width=2em]
  {
     1 & 2\\};
  \path[->,-stealth]
    (m-1-1) edge [bend left] node [above] {$a$} (m-1-2)
    (m-1-1) edge [bend right] node [below] {$b$} (m-1-2);
\end{tikzpicture}
\end{center}

It is well-known (see e.g. \cite[Theorem~4.3.2]{Benson1998RepresentationsCohomologyI} or \cite{Burgermeister1986ClassificationRepresentationsDoubleFleche}) that the indecomposable preprojective and preinjective $k$-representations of~$Q$ are given by 
\begin{center}
\begin{tikzpicture}
\matrix (m) [matrix of math nodes,row sep=1em,column sep=3em]
  {
     P_i: & k^{i} & k^{i+1}, && Q_i: & k^{i+1} & k^{i},\\};
  \path[->,-stealth]
    (m-1-2) edge [bend left] node [above] {$\begin{bsmallmatrix} \id \\ 0\end{bsmallmatrix}$} (m-1-3)
			edge [bend right] node [below] {$\begin{bsmallmatrix} 0 \\ \id\end{bsmallmatrix}$} (m-1-3)
    (m-1-6) edge [bend left] node [above] {$\begin{bsmallmatrix} \id & 0\end{bsmallmatrix}$} (m-1-7)
			edge [bend right] node [below] {$\begin{bsmallmatrix} 0 & \id\end{bsmallmatrix}$} (m-1-7);
\end{tikzpicture}
\end{center}
respectively, while the indecomposable regular representations are given by

\begin{center}
\begin{tikzpicture}
\matrix (m) [matrix of math nodes,row sep=1em,column sep=3em,minimum width=2em]
  {
     R_n(\phi,\psi) \colon & k^n & k^n,\\};
  \path[->,-stealth]
    (m-1-2) edge [bend left] node [above] {$\phi$} (m-1-3)
			edge [bend right] node [below] {$\psi$} (m-1-3);
\end{tikzpicture}
\end{center}
where either $\phi$ is the identity and $\psi$ is given by the companion matrix of a power of a monic irreducible polynomial over~$k$, or $\psi$ is the identity and $\phi$ is given by the companion matrix of a polynomial of the form $\lambda^m$.

We will show that the preprojective, the regular and the preinjective component are each hyperfinite families to conclude the amenability of $\mods kQ$. We will give an argument for the indecomposable objects in each component and then apply Proposition~\ref{prop:AdditiveClosureStaysHyperfinite} to extend the result.

We start with the preprojectives, and let $\varepsilon > 0$. Set $K_\varepsilon := \ceil*{\frac{1}{2\varepsilon}} + 1$ and $L_\varepsilon = \frac{1}{\varepsilon} + 3$. Let $X = P_i$ be some indecomposable preprojective. If $\dim X \leq L_\varepsilon$, there is nothing to show. We may thus assume that $\dim X > L_\varepsilon$, implying $i \geq K_\varepsilon$, and write $i = j \cdot K_\varepsilon + r,$ where $0 \leq r < K_\varepsilon$. Now consider the standard basis $\{e_1, e_2, \dots e_i\}$ of~$k^i$. Let $U$ be the submodule of~$X$ generated by the subset
\begin{align*} & \{e_1, \dots, e_{K_\varepsilon-1}\} \cup \{e_{K_\varepsilon+1}, \dots, e_{2 K_\varepsilon - 1}\} \cup \dots \\ 
\cup & \{e_{(j-1) K_\varepsilon + 1}, \dots, e_{j K_\varepsilon - 1}\} \cup \{ e_{j K_\varepsilon + 1}, \dots, e_{i} \},\end{align*}
dropping  every $K_\varepsilon$-th basis vector at the source. Then $U$ decomposes into $j$ direct summands of type $P_{K_\varepsilon-1}$ and a smaller rest term. All summands will thus be of $k$-dimension smaller than~$2 (K_\varepsilon - 1) + 1 < L_\varepsilon$. Moreover,
\begin{align*}
\dim U &= \dim X - j  = \dim X - \frac{i - r}{K_\varepsilon} = \dim X - \frac{\dim X -1}{2 K_\varepsilon} + \frac{r}{K_\varepsilon} \\
& \geq \dim X - \varepsilon (\dim X - 1) > (1-\varepsilon) \dim X.
\end{align*}
This shows that the family of indecomposable preprojective modules $\cP(kQ)$ is hyperfinite.

If $X = R_n(\phi,\psi)$ is an indecomposable regular module, we may consider the submodule $Y$ generated by the basis vectors $\{e_1, \dots ,e_{n-1}\}$ of the vector space at vertex $1$. Note that we assume that $\psi$ corresponds to the Frobenius companion matrix of a monic polynomial. Then $Y \isom P_{n-1}$, so by the above it belongs to the hyperfinite family $\cP(kQ)$. We have that $\dim Y  = \dim X - 1$. Thus, Proposition~\ref{prop:ExtendingHFfromSubmodulesOfBoundedCodimension} implies the hyperfiniteness of the indecomposable regular modules.

We are left to deal with the preinjective case. By Lemma~\ref{lemma:DescentOnDefectOfPIs}, for each indecomposable preinjective $X$, we can find a submodule $Y := \ker \theta$ of strictly smaller defect. Moreover, if $Y$ had a preinjective summand $Z$, it must have defect $\del(Z) < \del(X)$.
In this situation, all indecomposable preinjective modules have defect $d=1$. Choose the hyperfinite family $\cN_0 = \cP(kQ) \cup \cR(kQ)$ of all preprojective and regular modules. For all preinjective indecomposables, the submodule $Y$ must be in $\add \cN_0$, since there are no non-zero preinjective modules $Z$ with defect $\del(Z) < 1$. This family is hyperfinite by the above. Moreover, the codimension of~$Y$ is bounded by the dimension of the indecomposable injectives, of which there are only two. Hence, we can use Proposition~\ref{prop:ExtendingHFfromSubmodulesOfBoundedCodimension} to prove the hyperfiniteness of the indecomposable preinjectives.

Now apply Proposition~\ref{prop:AdditiveClosureStaysHyperfinite} to $\cP(kQ) \cup \cR(kQ) \cup \cQ(kQ)$ to see that $\mods kQ$ is hyperfinite, and thus $kQ$~is amenable.
\end{proof}
\end{theorem}

\section{Amenability of Extended Dynkin quivers} \label{section:ExtendedDynkinAmenable}

In order to prove our main result, another result is needed.

\begin{proposition} \label{prop:PerpendicularCategoryHF} \label{prop:PerpCatIsQuiverRep}
Let $Q$ be a finite acyclic quiver. 
\begin{enumerate}
\item If $T \in \mods kQ$ is an exceptional module without preprojective summands, $T^\perp$ is equivalent to $kQ'$ for some quiver $Q'$.
\item Assume $Q$ is of tubular type $(p,q,r)$, where $p > 1$, and that the extended Dynkin quiver of type $(p-1,q,r)$ is amenable. If $T$ is an inhomogeneous simple regular module belonging to a tube of rank $p$ in~$\Gamma_{kQ}$, then $T^{\perp}$ is hyperfinite.
\end{enumerate}
\begin{proof}
By \cite[Proposition~1.1]{GeigleLenzing1991PerpendicularCategoriesApplications}, in either case $T^{\perp}$ is an exact abelian subcategory of $\mods kQ$ closed under the formation of kernels, cokernels and extensions. What is more, \cite[Theorem~4.16]{GeigleLenzing1991PerpendicularCategoriesApplications} yields that $T^{\perp} = \mods \Lambda$ for some finite dimensional hereditary algebra $\Lambda$, along with a homological epimorphism $\varphi \colon kQ \to \Lambda$, which induces a functor $\varphi_* \colon \mods \Lambda \to \mods kQ$.
By Morita equivalence, we may assume that $\Lambda$ is basic.

Now, if $S$ is any simple $\Lambda$-module, then $S \isom P/\rad(P)$, where $P$ is a principal indecomposable $\Lambda$-module. By \cite[Theorem~4.4]{GeigleLenzing1991PerpendicularCategoriesApplications}, the natural maps $\End_{\Lambda}(P) \to \End_{kQ}(\varphi_* P)$ and $\Ext{\Lambda}{1}(P,P) \to \Ext{kQ}{1}(\varphi_* P,\varphi_* P)$ are isomorphisms, so $\varphi_*$ maps exceptional modules to exceptional modules. It follows from \cite[Corollary~1]{Ringel1994BraidGroupActionExceptionalSequences} that 
$\End_{kQ}(\varphi_* P) \isom \End_{kQ}(E)$ for some simple $kQ$-module $E$. But the simple $kQ$-modules all have trivial endomorphism ring~$k$.
Next, note that by \cite[Theorem~17.12]{AndersonFuller1992RingsCategoriesModules}
\[\End_{\Lambda}(S) \isom \End_{\Lambda}(P) / \J(\End_{\Lambda}(P)),\]
where $\J$ denotes the Jacobson radical. Recall that $\End_{\Lambda}(P) \isom \End_{kQ}(\varphi_* P) \isom k$,
thus $\J(\End_{\Lambda}(P)) = 0$. Hence it follows that $\End_{\Lambda}(S) \isom k$.
At large,
\[\End(\Lambda) / \J(\End(\Lambda)) \isom k \times \dots \times k.\]
Finally, \cite[Proposition~III.1.13]{AuslanderReitenSmaloe1995RepresentationTheoryArtinAlgebras} shows that $\Lambda$ is isomorphic to~$kQ'$ for some quiver~$Q'$.

It remains to prove the additional statements of (2). The proof of \cite[Theorem~13]{Schofield1986UniversalLocalisation} implies that $\Lambda$ is tame and hence the path algebra of an extended Dynkin quiver, and has tubular type $(p-1,q,r)$. By the hypothesis, it is amenable.

Now, if $F \colon \mods kQ' \to \mods \Lambda \to T^\perp$ is an equivalence, the simples $S(i)$ of~$kQ'$ get mapped to certain modules $B_i$ in $\mods kQ$.
The $k$-dimension of any module $M$ over a path algebra is determined by the length of any composition series. Such a series for $M$ in~$kQ'$ gets mapped to a composition series in the perpendicualar category, and thus a series in $\mods kQ$, such that the factor modules are isomorphic to some $B_i$. Letting $K' := \max\{\dim B_i\}$, we thus know that \[\dim_{k} F(M)_{kQ} \leq K' \dim_k M_{kQ'}.\]
On the other hand, if $F(M) \in T^\perp$, any submodule of~$F(M)$ in~$T^\perp$ is also a submodule in $\mods kQ$, so a composition series of~$F(M)$ in~$\mods kQ$ is at least as long as one in $T^\perp$. Thus, \[\dim_k M_{kQ'} \leq \dim_k F(M)_{kQ},\] using the fact that the length of~$M$ in $\mods kQ'$ equals the length of~$F(M)$ considered as an object of~$T^\perp$.
Hence by Proposition~\ref{prop:HFPreservingFunctors}, we have that each $T^{\perp}$ is a hyperfinite family.
\end{proof}
\end{proposition}

\begin{remark}
The above proof shows a slight improvement of \cite[Thm.~10.1(3)]{GeigleLenzing1991PerpendicularCategoriesApplications}, by removing the condition on~$k$ to be algebraically closed.
\end{remark}

\begin{theorem} \label{thm:ExtendedDynkinAmenable}
Let~$Q$ be an acyclic quiver of extended Dynkin type $\widetilde{A}_n$, $\widetilde{D}_n$, $\widetilde{E}_6$, $\widetilde{E}_7$ or~$\widetilde{E}_8$. Let~$k$ be any field. Then the path algebra~$kQ$ of~$Q$ is of amenable representation type.
\begin{proof}
Recall the tubular types and minimal radical vectors~$h_Q$ of the extended Dynkin diagrams in Table~\ref{tab:TubularTypes}, e.g. from \cite[p.~335]{Ringel1979InfiniteDimensionalRepresentationsFDHereditaryAlgebras}. Note that $\widetilde{A}_{p,q}$~is the quiver of type~$\widetilde{A}_{n}$ with $p+q = n+1$~vertices, where there are $p$~arrows in clockwise and $q$~arrows in anti-clockwise orientation.

\begin{table} 
\begin{center}
\begin{tabular}{l | c c}
$Q$ & $(m_i)$ & $h_Q$ \\ \hline\noalign{\smallskip}
$\widetilde{A}_{p,q}$ & $(p,q)$ & $\begin{smallmatrix} & 1 \dots 1 \\ 1 && 1 \\ & 1 \dots 1 \end{smallmatrix}$ \\[10pt]
$\widetilde{D}_n$ & $(2,2,n-2)$ & $\begin{smallmatrix} 1 &&&& 1\\ & 2 & \dots & 2 \\ 1 &&&& 1 \end{smallmatrix}$ \\[10pt]
$\widetilde{E}_6$ & $(2,3,3)$ & $\begin{smallmatrix} && 1\\ && 2 \\ 1 & 2 & 3 & 2 & 1 \end{smallmatrix}$ \\[10pt]
$\widetilde{E}_7$ & $(2,3,4)$ & $\begin{smallmatrix} &&& 2 \\ 1 & 2 & 3 & 4 & 3 & 2 & 1 \end{smallmatrix}$ \\[10pt]
$\widetilde{E}_8$ & $(2,3,5)$ & $\begin{smallmatrix} && 3\\ 2 & 4 & 6 & 5 & 4 & 3 & 2 & 1 \end{smallmatrix}$
\end{tabular}
\end{center}
\caption{Tubular types and minimal radical vector of the acyclic extended Dynking diagrams.} \label{tab:TubularTypes}
\end{table}
 
We will prove the claim by induction on~$n$ for the case of the acyclic~$\widetilde{A}_n$ and for~$\widetilde{D}_n$, and use the case of~$\widetilde{D}_5$ to prove it for the $\widetilde{E}$-family, stepping from $5$ to~$6$ to~$7$.

\proofstep{Case~$\widetilde{A}_n$:}
For~$\widetilde{A}_1$, the only acyclic orientation is the $2$-Kronecker quiver, which has been shown to be of amenable type in Theorem~\ref{thm:2KroneckerAmenable}.

Now assume the case of all acyclic quivers~$\widetilde{A}_n$ for~$n \geq 1$ has already been proven. Let~$\widetilde{A}_{p,q}$ be of type $\widetilde{A}_{n+1}$. Then $p+q = n+2 \geq 3$. We may thus assume that $p \geq 2$, and choose a tube $\T$ of rank $m := p$, and denote the isoclasses of simple regular modules in this tube by $T_1, \dots, T_m$. Since all indecomposable preprojective~$kQ$-modules $X$ have defect $\del(X) = -1$, Lemma~\ref{lemma:PPofDefectSmallerThanTubeRank} implies that every indecomposable preprojective is contained in the perpendicular category $T_i^{\perp}$ for some $1 \leq i \leq m$. 
By Proposition~\ref{prop:PerpendicularCategoryHF}, each $T_i^{\perp}$ is hyperfinite. This shows that the preprojectives form a hyperfinite family, using Proposition~\ref{prop:AdditiveClosureStaysHyperfinite}.

Next, we consider the regular modules. Indecomposable regular modules in a tube other than $\T$ will be contained in $T_1^{\perp}$ by \cite[3.1.(3')]{Ringel1984TameAlgebrasIntegralQuadraticForms}. 
By Lemma~\ref{lemma:RegularsInTubeRankGreaterOne}, any regular indecomposable in $\T$ either is contained in the perpendicular category of some simple regular in $\T$ or has a submodule of bounded codimension, which is in the perpendicular category of some simple regular in $\T$. But by the above argument, the perpendicular categories are hyperfinite. In the latter case, we can therefore apply Proposition~\ref{prop:ExtendingHFfromSubmodulesOfBoundedCodimension} to show the hyperfiniteness of these indecomposable regular modules.

For the preinjective modules, we may apply the argument used for the preinjectives in the proof of Thereom~\ref{thm:2KroneckerAmenable}. This completes the induction step.

\proofstep{Case~$\widetilde{D}_n$:}
For the case of~$\widetilde{D}_4$, choose a tube $\T$ of rank $2$, and denote the simple regular modules in $\T$ by $S$ and $T$.
In this case, as the extended Dynkin quiver of tubular type $(1,2,2)$ is one of type~$\widetilde{A}_{2,2}$, which is known by the above to be amenable, Proposition~\ref{prop:PerpendicularCategoryHF} implies that $S^{\perp}$ and $T^{\perp}$ are hyperfinite.

All preprojective modules $X$ of defect $\del(X) = -1$ are in $S^{\perp}$ or $T^{\perp}$ by Lemma~\ref{lemma:PPofDefectSmallerThanTubeRank}.
Using Lemma~\ref{lemma:PPofDefectSmallerThanTwiceTubeRank}, we can find a submodule $Y$ for all but finitely many indecomposable preprojectives $X$ of defect $\del(X) = -2$, which are not themselves in $S^{\perp}$ or $T^{\perp}$. Since the dimension vector of the simple regular in $\T$ is bounded, the conditions of Propositon~\ref{prop:ExtendingHFfromSubmodulesOfBoundedCodimension} are satisfied for all but finitely many indecomposable preprojectives of defect~$-2$. This shows that the preprojectives form a hyperfinite family.

Moreover, the regular modules are hyperfinite: If they are in a tube other than $\T$, they will be contained in $S^{\perp}$ by \cite[3.1.(3')]{Ringel1984TameAlgebrasIntegralQuadraticForms}. Choosing a second non-homogeneous tube $\T'$ and a simple regular $U \in \T'$, we know that $\T \subset U^\perp$, which is also hyperfinite.

We are left to deal with the preinjective modules. By Lemma~\ref{lemma:DescentOnDefectOfPIs}, for each indecomposable preinjective $X$, we can find a submodule $Y := \ker \theta$ of strictly smaller defect. Moreover, if $Y$ had a preinjective summand $Z$, it must have defect $\del(Z) < \del(X)$. We do an induction on the defect~$d$. If $d=1$, then we can choose the hyperfinite family $\cN_0$ of all preprojective and regular modules. For all preinjective indecomposables of defect~$d$, the submodule $Y$ must be in $\add \cN_0$, since there are no preinjective modules with defect $\del(Z) < 1$. This family is hyperfinite by Proposition~\ref{prop:AdditiveClosureStaysHyperfinite} and the above results. Moreover, the codimension of~$Y$ is bounded by the dimension of the indecomposable injectives, of which there  are only finitely many. Hence, we can use Proposition~\ref{prop:ExtendingHFfromSubmodulesOfBoundedCodimension} to prove the hyperfiniteness of the indecomposable preinjectives of defect~one.
We recursively define \[\cN_d := \cN_{d-1} \cup \{\text{indecomposable preinjectives of defect~$d$}\}.\] Note that the base case implies that $\cN_1$ is hyperfinite.
For the induction, note that Lemma~\ref{lemma:DescentOnDefectOfPIs} also yields a submodule in $\add \cN_d$ for every indecomposable preinjective of defect~$d+1$ of bounded codimension. Assuming the hyperfiniteness of~$\cN_d$, Proposition~\ref{prop:ExtendingHFfromSubmodulesOfBoundedCodimension} yields that $\cN_{d+1}$ is hyperfinite.
This proves the claim for~$\widetilde{D}_4$, using Proposition~\ref{prop:AdditiveClosureStaysHyperfinite}.

Now assume the case of~$\widetilde{D}_{n}$ has already been dealt with for $n \geq 4$. To prove the amenability of~$\widetilde{D}_{n+1}$, choose $\T$ to be the unique tube of maximal rank. Similarly to the base case, $S^{\perp}$ is of amenable type for $S \in \T$, since the tubular type $\left(2,2,(n+1)-2 - 1\right)$ belongs to~$\widetilde{D}_{n}$.
By inspection of Table~\ref{tab:TubularTypes}, we see that the indecomposable preprojectives have negative defect~one or~two. Hence, they are in a hyperfinite family by Lemma~\ref{lemma:PPofDefectSmallerThanTubeRank}. The regular indecomposables are hyperfinite by an argument similar to that of the base case. To deal with the indecomposable preinjectives, we may apply the same induction as in the base case.

\proofstep{Case~$\widetilde{E}_n$}
We proceed with~$\widetilde{E}_n$ for $n=6,7,8$. Assume the path algebra of tubular type $(2,3,n-4)$ has already been shown to be of amenable type. By choosing $\T$ to be a tube of maximal rank $m=n-3$, we find simple regular modules $S$ such that $S^{\perp}$ is hyperfinite, for the argument of Proposition~\ref{prop:PerpendicularCategoryHF} reduces it to $(2,3,m-1)$. Inspection tells us that any indecomposable preprojective module will have negative defect less than~$2m$. Thus, we can use Lemma~\ref{lemma:PPofDefectSmallerThanTwiceTubeRank} -- if needed in connection with Proposition~\ref{prop:ExtendingHFfromSubmodulesOfBoundedCodimension} -- to show that all but finitely many, and hence all preprojective indecomposables form a hyperfinite family. For the indecomposable regular and preinjective modules, use the same arguments as for~$\widetilde{D}_n$.
\end{proof}
\end{theorem}

\section{Non-amenability of wild path algebras} \label{section:WildQuivers}
We say that a quiver $Q$ is wild if it strictly contains an extended Dynkin quiver, i.e. has indefinite quadratic form.

\begin{prop} \label{prop:NonAmenabilityFromSubquiver}
Let $Q$ be a quiver. If $Q$ has a subquiver $Q'$ such that $\mods kQ'$ is not amenable, then $\mods kQ$ is not of amenable representation type.
\begin{proof}
Let $F \colon \mods kQ' \to  \mods kQ$ be the embedding mapping a representation $M' = ((M'_i),(M'_a))$ of $Q'$ to the representation $M = ((M_i),(M_a))$ of $Q$ given by
\[M_i = \begin{cases} M'_i, & i \in Q'_0,\\ 0, & \text{else,} \end{cases}
\qquad	
M_a = \begin{cases} M'_a, & a \in Q'_1,\\ 0, & \text{else.}\end{cases}\]
Since $\mods kQ'$ is not of amenable type, there exists  a non-hyperfinite family of modules $\{M_n\}_{n=1}^{\infty} \subseteq \mods kQ'$. Assume that $\{FM_n\}_{n \in \N}$ is hyperfinite, for otherwise  we have found a non-hyperfinite sequence exhibiting the non-amenability of $\mods kQ$.
Note that any submodule $N$ of some $FM_n$ is given by subspaces $N_i \subseteq (FM_n)_i$ for each $i \in Q_0$ and linear maps $N_a$ for each $a \in Q_1$ such that $N_a(N_{s(a)}) \subseteq N_{t(a)}$.
So $N' := ((N_i)_{i \in Q'_0},(N_a)_{a \in Q'_1})$ is a subrepresentation of $M_n$.
Moreover, $\dim_k FM_n = \dim_k M_n$ and $\dim_k N' = \dim_k N$. 
Also, if $S$ is a direct summand of $N$, then each $S_i$ is a direct summand of $N_i$. This fact along with $S_i = 0$ for all $i \in Q_0 \setminus Q'_0$ implies that $S$ also yields a direct summand $S'$ of $N'$ of dimension $\dim_k S' =\dim_k S$.
Alltogether, this implies that $\{M_n\}_{n \in \N}$ is hyperfinite, a contradiction.
\end{proof}
\end{prop}

\begin{prop} \label{prop:NonHFPreservingFunctors}
Let $k$ be a field and $A,B$ be two finite dimensional $k$-algebras. Let $\{M_n \colon n \in \N\} \subseteq \mods A$ be a non-hyperfinite family of modules.
If there exist $K_1, K_2 > 0$ and additive functors $F \colon \mods A \to \mods B$ and $G \colon \mods B \to \mods A$ such that 
\begin{itemize}
\item $G F (M_n) \isom M_n$ for all $n \in \N$,
\item $G$ is left exact,
\item $K_1 \dim F (M_n) \leq \dim G F (M_n)$ for all $n \in \N$,
\item $\dim G(X) \leq K_2 \dim X$ for all $X \in \mods B$,
\end{itemize}
then $\{F(M_n) \colon n \in \N\}$ is a non-hyperfinite family, too.
\begin{proof}
Consider the family $\{F(M_n) \colon n \in \N\}$ in $\mods B$. Assume that it is hyperfinite. By Proposition~\ref{prop:HFPreservingFunctors} we know that $\{G F(M_n) \colon n \in \N\} = \{M_n \colon n \in \N\}$ is hyperfinite, as the ensuing remark shows that (3) is sufficient to apply the proposition. This yields a contradiction, so $\{F(M_n) \colon n \in \N\}$ cannot be hyperfinite, exhibiting the non-amenability of $\mods B$.
\end{proof}
\end{prop}

\begin{theorem} \label{thm:WildQuiverIsNotAmenable}
Let $k$ be a finite field and $Q$ a wild acyclic quiver. Then $\mods kQ$ is not of amenable representation type.
\begin{proof} 
If $Q$ contains the wild 3-Kronecker quiver as a subquiver, Proposition~\ref{prop:NonAmenabilityFromSubquiver} along with \cite[Theorem~6]{Elek2017InfiniteDimensionalRepresentationsAmenabilty} implies that $\mods kQ$ is not of amenable type.

If $Q$ is a wild quiver but does not contain a wild Kronecker quiver, it must have $n \geq 3$ vertices, hence $\mods kQ$ has at least three isoclasses of simple modules. By \cite[Theorem~2.1]{Baer1989NoteWildQuiverAlgebrasTiltingModules}, which also holds for arbitrary fields, there exists a regular indecomposable module $S$ without self-extensions. Thus we can apply Proposition~\ref{prop:PerpCatIsQuiverRep} to see that $S^\perp$ is isomorphic to $\mods kQ'$ for some quiver $Q'$. Recall from the proof thereof that there is a corresponding homological epimorphism $\varphi \colon kQ \to kQ'$ and the induced functor $F = \varphi_* \colon \mods kQ' \to \mods kQ$ is fully faithful and exact.
Indeed, $Q'$~has $n-1$ vertices and \cite[Theorem~4.1]{Baer1989NoteWildQuiverAlgebrasTiltingModules} shows that $kQ'$ is a wild quiver algebra.

Now, if we have some equivalence $\mods kQ' \iso S^\perp$, the simples $S(i)$ of $kQ'$ are mapped to certain objects $B_i$, considered as modules in $\mods kQ$.
The $k$-dimension of any module $M$ over a path algebra is determined by the length of any composition series. Such a series for some $M'$ in $kQ'$ is mapped to a composition series in the perpendicular category, and thus a series in $\mods kQ$, such that the factor modules are isomorphic to some $B_i$. This shows that \[\dim_{k} F(M')_{kQ} \leq \max_{i=1,\dots,n-1}\{\dim_k B_i\} \dim_k {M'}_{kQ'}.\]

Moreover, by \cite[Proposition~1.1]{GeigleLenzing1991PerpendicularCategoriesApplications}, $S^\perp$ is closed under kernels, cokernels and extensions. Assuming that $Z$ is a relative projective generator of $\mods kQ'$, \cite[Proposition~A.1]{KrauseHubery2016CategorificationNonCrossingPartitions} shows that there exists a right adjoint to the inclusion $F$, which we denote by $G \colon \mods kQ \to \mods kQ'$. Note that $G$ is left exact.
Indeed, for $M \in \mods kQ$, $G(M)$ is given as a factor module of a right $\add(Z)$-approximation $Z_M$ of $M$. Since we may assume that $Z_M \isom Z \otimes_{\End(Z)} \Hom(Z,M)$, this implies that
\[\dim_k G(M) \leq \dim_k Z_M \leq (\dim_k Z)^2 \dim_k M.\]
Since $F$ is fully faithful, we have that  $X \iso GF(X)$ for all $X \in \mods kQ'$.

To conclude the proof, just note that we have prepared a descent argument leading to a wild quiver with two vertices, which has to include a wild Kronecker quiver $\Theta(m)$ as a subquiver. 
Given any non-hyperfinite sequence $\{M_n\}$ in $\mods k\Theta(m)$, we choose $Z$ as above and let $K_1 = \max\{\dim_k B_i\}^{-1}$ and $K_2 = (\dim_k Z)^2.$ We can now choose all $F_n$ as $F$ and all $G_n$ as $G$ and have fulfilled the conditions of Proposition~\ref{prop:NonHFPreservingFunctors}, which we apply to show the non-hyperfiniteness in each step, until we reach $\mods kQ$.
\end{proof}
\end{theorem}

\begin{remark}
The above theorem can also be proved more directly, by supplying a tangible embedding for each minimal wild quiver $Q$
(see \cite[Section~4]{Kerner1988PreprojectiveComponentsWildTiltedAlgebras} for a definition and a list). In each case, we need to exhibit exceptional objects $X,Y \in \mods kQ$ such that $(X,Y)$ is an orthogonal exceptional pair, i.e. \[\Hom(Y,X) = \Hom(X,Y) = 0 = \Ext{}{1}(Y,X),\] and such that $m := \dim_k \Ext{}{1}(X,Y) \geq 3$.

Then, by \cite[1.5 Lemma]{Ringel1976RepresentationsKSpeciesBimodules}, there is a full exact embedding $F \colon \mods k\Theta(m) \to \mods kQ$.
This will map the simple representations of $\Theta(m)$ to $X$ and $Y$ respectively. Now, if $M$ is a module for $k\Theta(m)$, any composition series will get mapped to a series in $\mods kQ$, such that the factor modules are isomorphic to either $X$ or $Y$. This shows that
\[\dim_k F(M) \leq \max\{\dim_k X,\dim_k Y\} \dim_k M. \]

Denoting the closure of the full subcategory of $\mods kQ$ containing $X$ and $Y$ under kernels, images, cokernels and extensions by $\C(X,Y)$, $F$ induces an equivalence $\mods k\Theta(m) \iso \C(X,Y)$, see for example \cite[Section~1]{Ringel1976RepresentationsKSpeciesBimodules} in connection with \cite[Corollary~1]{Ringel1994BraidGroupActionExceptionalSequences}.
Assuming that $Z$ is a relative projective generator of $\C(X,Y)$, \cite[Proposition~A.1]{KrauseHubery2016CategorificationNonCrossingPartitions} shows that there exists a right adjoint to the inclusion, which we denote by $G \colon \mods kQ \to \C(X,Y)$. Moreover, if $M \in \mods kQ$, $G(M)$ is given as a factor module of a right $\add(Z)$-approximation $Z_M$ of $M$. Since we may assume that $Z_M \isom Z \otimes_{\End(Z)} \Hom(Z,M)$, this implies that
\[\dim_k G(M) \leq \dim_k Z_M \leq (\dim_k Z)^2 \dim_k M.\] Moreover, $G$ is left exact and we have $GF(M) \isom M$ for all $M \in \mods k\Theta(m)$.

To conclude the proof in this case, we may chose a non-hyperfinite sequence $\{M_n\}$ in $\mods k\Theta(m)$, guaranteed to exist by \cite[Theorem~6]{Elek2017InfiniteDimensionalRepresentationsAmenabilty}. Let $X,Y,Z$ as above and set $K_1 = \max\{\dim_k X,\dim_k Y\}^{-1}$ and $K_2 = (\dim_k Z)^2.$ We can now choose all $F_n$ as $F$ and all $G_n$ as $G$ and have fulfilled the conditions of Proposition~\ref{prop:NonHFPreservingFunctors}, which we apply to show the non-hyperfiniteness of $\mods kQ$. 
\end{remark}

\begin{examples}
\begin{tikzpicture}[baseline={($(current bounding box.north)-(0,1.6ex)$)}]
\matrix (a) [matrix of math nodes, row sep=2em, column sep=1em,label={above left:{(1)}}]
  {
	&&0\\
	1 & 2 & 3  & 4 & 5\\
  };
  \path[<-]
	(a-1-3) edge (a-2-1)
			edge (a-2-2)
			edge (a-2-3)
			edge (a-2-4)
			edge (a-2-5);
	
\matrix (d) [matrix of math nodes, text height=1.5ex, text depth=0.25ex, row sep=1.5em, column sep=1.5em, right=1cm of a,label={above left:{(2)}}]
  {
	  &   & 6\\
	  &   & 5\\
	1 & 2 & 3 & 4 & 0 & \infty\\
  };
  \path[->]
	(d-3-1) edge (d-3-2)
	(d-3-2)	edge (d-3-3)
	(d-3-4)	edge (d-3-3)
	(d-3-5) edge (d-3-4)
	(d-3-6) edge (d-3-5)
	(d-1-3) edge (d-2-3)
	(d-2-3) edge (d-3-3);
\end{tikzpicture}

\begin{enumerate}
\item Let $Q$ be the five subspace quiver $\bS(5)$. Choose the modules $X = S(5)$, the simple for vertex $5$ and $Y = \tau^{-1}P(5)$. Then \[\Hom(X,Y) = 0 = \Hom(Y,X),\] since $X$ and $Y$ have disjoint support. Also, $\Ext{}{1}(Y,X)=0$ since $X$ is injective. On the other hand, $\Ext{}{1}(X,Y) \isom k^3$ for dimension reasons.
\item Let $Q$ be the one-point extension at an extending vertex $0$ of $\tilde{E}_6$ oriented towards the center by a source $\infty$. Choose $X=S(\infty)$ and $Y$ as the representation induced by $\tau^{-6} P(1)$ of the underlying $\tilde{E}_6$.
\end{enumerate}

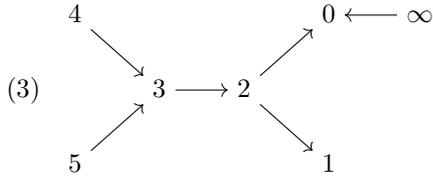
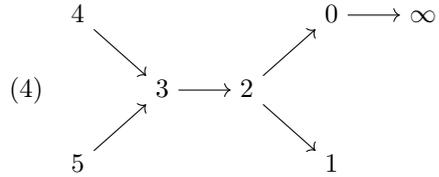
\centering
\begin{tikzpicture}
\matrix (b) [matrix of math nodes, text height=1.5ex, text depth=0.25ex, row sep=1.41em, column sep=2em, label={left:{(3)}}]
  {
	4 &   &   & 0 & \infty\\
	  & 3 & 2\\
	5 &   &   & 1\\
  };
  \path[->]
	(b-1-1) edge (b-2-2)
	(b-3-1)	edge (b-2-2)
	(b-2-2)	edge (b-2-3)
	(b-2-3) edge (b-1-4)
			edge (b-3-4)
	(b-1-5) edge (b-1-4);	
	
\matrix (c) [matrix of math nodes, text height=1.5ex, text depth=0.25ex, row sep=1.41em, column sep=2em, right=2cm of b,label={left:{(4)}}]
  {
	4 &   &   & 0 & \infty\\
	  & 3 & 2\\
	5 &   &   & 1\\
  };
  \path[->]
	(c-1-1) edge (c-2-2)
	(c-3-1)	edge (c-2-2)
	(c-2-2)	edge (c-2-3)
	(c-2-3) edge (c-1-4)
			edge (c-3-4)
	(c-1-4) edge (c-1-5);
\end{tikzpicture}
\begin{enumerate}
\addtocounter{enumi}{2}
\item Let $Q$ be the one-point extension at an extending vertex $0$ of a linearly oriented $\tilde{D}_5$ by a source $\infty$. Choose $X=S(\infty)$ and $Y$ as the representation induced by $\tau^{-3}P(3)$ of the underlying $\tilde{D}_5$.
\item Let $Q$ be the one-point extension at an extending vertex $0$ of a linearly oriented $\tilde{D}_5$ by a sink $\infty$. Choose $Y=S(\infty)$ and $X$ as the representation induced by $\tau^{-3}P(3)$ of the underlying $\tilde{D}_5$.
\end{enumerate}
\end{examples}

\clearpage
\begin{acknowledgement}
This note is based on work done during the author’s doctorate studies at Bielefeld University. The author would like to thank his supervisor Professor W. Crawley-Boevey for his advice and guidance and A. Hubery for helpful discussions.
\end{acknowledgement}

{
 \printbibliography
}
\end{document}